\documentclass[a4paper,12pt,leqno]{article}
\usepackage{latexsym}
\usepackage[all]{xy}

\usepackage{amssymb} 
\usepackage{amsmath} 

\usepackage{tikz}
\usetikzlibrary{arrows,cd}  
\usepackage{amsthm}

\def\P{{\mathbf{P}}}
\def\Z{{\mathbb{Z}}}
\def\K{{\mathbb{K}}}
\def\CC{{\mathbb{C}}}
\def\R{{\mathbb{R}}}
\def\A{{\mathcal{A}}}
\def\B{{\mathcal{B}}}

\DeclareMathOperator{\coker}{coker}
\DeclareMathOperator{\Der}{Der}

\DeclareMathOperator{\res}{res}

\DeclareMathOperator{\POexp}{POexp}

\DeclareMathOperator{\hilb}{Hilb}  

\numberwithin{equation}{section}

\newcommand{\owari}{\hfill$\square$}

\newtheorem{theorem}{Theorem}[section]
\newtheorem{prop}[theorem]{Proposition}
\newtheorem{cor}[theorem]{Corollary}
\newtheorem{lemma}[theorem]{Lemma}
\newtheorem{define}[theorem]{Definition}

\newtheorem{problem}[theorem]{Problem}
\newtheorem{conj}[theorem]{Conjecture}
\theoremstyle{remark}

\newtheorem{example}[theorem]{Example}

\title{Cokernels of the Euler restriction map of logarithmic derivation modules}
\author{Takuro Abe and Hiraku Kawanoue}
\date{\today}

\begin{document}

\maketitle

\begin{abstract}There are two restriction maps of the logarithmic modules of plane arrangements in a three dimensional vector space. One is the Euler restriction and the other is the Ziegler restriction. The dimension of the cokernel of the Ziegler restriction map of logarithmic derivation modules has been well-studied for the freeness of 
hyperplane arrangements after 
Yoshinaga's celebrated criterion for freeness, which connects the second Betti number and the splitting type (exponents). 
However, though the Euler restriction has a longer history than the Ziegler restriction, the cokernel and its dimension of the Euler restriction have not been studied at all. 

The aim of this article is to study the cokernel and dimension of the Euler restriction maps in terms of combinatorics, more explicitly, the characteristic polynomial. We give an upper bound of that cokernel, and show the formula for that if the arrangement is free. 
  \end{abstract}

\section{Main theorem}

Let $\K$ be a field and $\A$ be a central arrangement in $\K^3=V$, i.e., a finite set of linear 
planes in $V$. Let $S=\mbox{Sym}(V^*)=\K[x_1,x_2,x_3]$ 
and let 
$\Der S:=\oplus_{i=1}^3 S \partial_{x_i}$.  
Fixing the defining linear form $\alpha_H \in V^*$ for each $
H \in \A$, we can define the \textbf{logarithmic derivation module} $D(\A)$ of $\A$ by 
$$
D(\A):=\{\theta \in \Der S \mid \theta(\alpha_H) \in S\alpha_H\ (\forall H \in \A)\}.
$$
Geometrically $D(\A)$ consists of polynomial vector fields of $V$ which are tangent to $\A$. 
$D(\A)$ is an $S$-graded module of rank three over $S$ and is known to be reflexive. Since $\dim V=3$, in this case 
$D(\A)$ is locally free, in particular, $\widetilde{D(\A)}$ is a vector bundle over 
$\mbox{Proj}(S)=\P^2$. When $\A \neq \emptyset$, then there is the 
\textbf{Euler derivation} $\theta_E=\sum_{i=1}^3 x_i \partial_{x_i}$ such that 
$$
D(\A)=S\theta_E \oplus D_H(\A),
$$
where 
$$
D_H(\A):=\{\theta \in D(\A) \mid 
\theta(\alpha_H)=0\}
$$
for each $H \in \A$. Thus $D_H(\A) \simeq 
D(\A)/S \theta_E \simeq D_L(\A)$ for any $H,L \in \A$. 
In general, $D(\A)$ is not free. Thus we say that $\A$ is \textbf{free} with $\exp(\A)=(1,a,b)$ if 
$$
D(\A)\simeq S[-1]\oplus S[-a]\oplus S[-b],
$$
or equivalently,
$$
D_H(\A)\simeq S[-a]\oplus S[-b].
$$
Here the degree $1$-part in the above corresponds to the free submodule $S\theta_E$ of $D(A)$. 

Let us define combinatorial invariants of $\A$. The \textbf{intersection lattice} $L(\A)$ of $\A$ 
is defined by 
$$
L(\A):=\{ \cap_{H \in \B} H \mid B \subset \A\}.
$$
$L(\A)$ knows how planes in $\A$ intersect, so we regard $L(\A)$ as combinatorial data of $\A$. We can define the \textbf{M\"{o}bius function} 
$\mu:L(\A) \rightarrow \Z$ defined by $\mu(V)=1$ and by $\mu(X)=-\sum_{ X \subsetneq Y \in L(\A)} 
\mu(Y)$ for $ X \neq V$. Then the \textbf{characteristic polynomial} $\chi(\A;t)$ of $\A$ is 
defined by 
$$
\chi(\A;t):=\sum_{X \in L(\A)} \mu(X)t^{\dim X}.
$$
When $\A \neq \emptyset$, it is known that $\chi(\A;t)$ is divisible by $t-1$. So 
$$
\chi_0(\A;t):=\chi(\A;t)/(t-1)=t^2-b_1^0t+b_2^0
$$
is defined, and called the \textbf{reduced characteristic polynomial} of $\A$. Both are the most 
important combinatorial invariants of $\A$. Note that it is easy to show by the definition of $\mu$ that 
$$
b_1^0=|\A|-1$$
and the coefficient of $t^2$ in $\chi(\A;t)$ is $-|\A|$. 

For each $H \in \A$, we can define two restrictions of $\A$ and $D(\A)$. The first one is the 
\textbf{Euler restriction} $\A^H:=\{L \cap H \mid L \in 
\A \setminus \{H\}\}$. Then in the same way we can define its logarithmic derivation module 
$D(\A^H)$, that is a module graded by $S/S\alpha_H$. Since $D(\A^H)$ is reflexive too and 
it contains $\theta_E$, we know that $\A^H$ is free with $\exp(\A^H)=(1,|\A^H|-1)$. $\A^H$ does not
contain any information of multiplicity, so we can define the other restriction called the \textbf{Ziegler restriction} $(\A^H,m^H)$, where $m^H:\A^H \rightarrow \Z_{>0}$ is defined by 
$$
m^H(X):=|\{L \in \A \setminus \{H\} \mid L \cap H=X\}|
$$
for each $X \in \A^H$. Then for $\overline{S}:=S/S\alpha_H$, we can define 
$$
D(\A^H,m^H):=\{\theta \in \Der \overline{S} \mid 
\theta(\alpha_X) \in \overline{S}\alpha_X^{m^H(X)}\ (\forall X \in \A^H)\}.
$$
Also for $\theta \in D(\A)$ the \textbf{Euler restriction map} $\rho^H:D(\A) \rightarrow D(\A^H)$ is defined by 
$$
\rho^H(\theta)(\overline{f}):=\overline{\theta(f)}
$$
for $f \in S$. Regarding $D_H(\A)$ as a submodule of $D(\A)$ by the above decomposition, the 
\textbf{Ziegler restriction map} $\pi^H:D_H(\A) \rightarrow 
D(\A^H,m^H)$ is defined by $\pi^H:=\rho^H|_{D_H(\A)}$. Then we have the \textbf{Euler exact sequence}
$$
0 \rightarrow D(\A \setminus \{H\}) \stackrel{\cdot \alpha_H}{\rightarrow}
D(\A) \stackrel{\rho^H}{\rightarrow}D(\A^H)$$
and
the \textbf{Ziegler exact sequence}
$$
0 \rightarrow D_H(\A) \stackrel{\cdot \alpha_H}{\rightarrow}
D_H(\A) \stackrel{\pi^H}{\rightarrow}D(\A^H,m^H).
$$
See \cite{OT} and \cite{Z} for details. 

Since $D(\A^H,m^H)$ is reflexive over $\overline{S}$ which is a coordinate ring of a two-dimensional vector space $H$, it is free. However, contrary to $D(\A^H)$, it is hard to determine 
$\exp(\A^H,m^H)$. This difficulty is similar to determine the splitting type of a rank two vector bundle over 
$\P^2$. These two restriction maps have played essential roles in the research of freeness. In particular, in \cite{Y2} Yoshinaga gave a criterion for freeness of $3$-arrangements by using Ziegler restriction. Namely, first, he computed the dimension of the cokernel of the Ziegler restriction map 
$
\pi^H:D_H(\A) \rightarrow D(\A^H,m^H).
$
It was known that $\pi^H$ is surjective if and only if $\A$ is free by Ziegler in \cite{Z}. Yoshinaga made this into numerical setup, i.e., if $\exp(\A^H,m^H)=(d_1,d_2)$, then 
$$
\dim \coker \pi^H=\chi_0(\A;0)-d_1d_2
$$
and it is zero if and only if $\A$ is free with $\exp(\A)=(1,d_1,d_2)$. 
Also, the fact that 
$\dim \coker \pi^H<\infty$ is frequently used and we obtain a lot of related results on freeness, 
combinatorics and geometry of hyperplane arrangements. Based on these developments, it is natural to ask the following first question in this article:

\begin{problem}
Can we compute the dimension of the cokernel of the 
Euler restriction map
$$
\rho^H:D(\A) \rightarrow D(\A^H)
$$
or not?
Is it of the finite dimensional at all?
\label{problem1}
\end{problem}

No one has considered this problem though $\rho^H$ was discovered and used far before $\pi^H$. This is the first problem we want to solve. On the other hand we consider a different problem from it, which seems rather combinatorial.

In \cite{A}, it was proved that 
$$
\chi_0(\A;|\A^H|-1) \ge 0
$$
for any $H \in \A$. Since $\chi_0(\A;t)$ is combinatorially determined and is a polynomial with integer coefficients, $\chi_0(\A;|\A^H|-1)$ is a combinatorially determined nonnegative integer. Then from the point of combinatorics, it is natural to ask whether there is some geometric or algebraic invariants that can be expressed by this integer. Let us give a name to it as follows:

\begin{define}
For $H \in \A$, the combinatorial invariant $\chi_0(\A;|\A^H|-1)$ is denoted by $LP(\A,H)$, called the \textbf{line-point invariant} or 
\textbf{LP-invariant}, which is a nonnegative integer, see Theorem \ref{roots}.
\label{LP}
\end{define}

So our next problem is as follows:

\begin{problem}
What is a combinatorial, algebraic or geometric meaning of $LP(\A,H)$?
\label{problem2}
\end{problem}

Our main result in this article surprisingly answers these two, completely different problems at once as follows:

\begin{theorem}
\begin{itemize}
    \item [(1)] It holds that 
    $$\dim \coker \rho^H \le \chi_0(\A;|\A^H|-1)=LP(\A;H),
$$ and
    \item[(2)] 
    $$\dim \coker \rho^H =LP(\A;H)$$
    if $\A$ is free with $\exp(\A)=(1,a,b)$. Explicitly,
    $$
    \hilb(\coker \rho^H;x)=x^{|\A^H|-1}
    \displaystyle \frac{1-x^{a-|\A^H|+1}}{1-x}
    \displaystyle \frac{1-x^{b-|\A^H|+1}}{1-x}.
    $$
\end{itemize}
\label{LPbound}
\end{theorem}

Thus the LP-invariant describes the upper bound of the dimension of the cokernel of the 
Euler restriction map, and it attains the upper bound when $\A$ is free. Hence Theorem \ref{LPbound} (1) gives a sharp upper bound. 

The organization of this article is as follows. After introducing basic results used for the proof of main results in \S2, we prove Theorem \ref{LPbound} in \S3. \S4 is devoted to understand Yoshinaga's criterion in terms of $\chi_0(\A;n)$ for some integer $n$, and give applications 
of Theorem \ref{LPbound}.
\medskip

\noindent
\textbf{Acknowledgements}.  The authors are grateful to Shizuo Kaji for his comment to Example 4.5. The first author 
is partially supported by JSPS KAKENHI Grant Numbers JP23K17298 and JP23K20788.
The second author is partially supported by JSPS KAKENHI Grant Number JP24K06656.

%



\section{Preliminaries}
In this section let us recall fundamental results on algebra and combinatorics for hyperplane arrangements to prove Theorem \ref{LPbound}. \textbf{For the rest of this paper, 
$\A$ is an arrangement in $\K^3$ unless otherwise specified}. First let us recall the famous Saito's freeness criterion.

\begin{theorem}[Saito's criterion, \cite{Sa}]
Let $H \in \A$, $Q(\A):=\prod_{H \in \A} \alpha_H$ and 
$Q(\A^H,m^H):=\prod_{X \in \A^H} \alpha_X^{m^H(X)}$.

\begin{itemize}

\item[(1)]
The homogeneous derivations $\theta_1:=\theta_E,\theta_2,\theta_3 \in D(\A)$ form a basis for $D(\A)$ if and only if 
$\det (\theta_i(x_j))=cQ(\A)$ up to nozero scalar $c \in \K$. In particular, 
$\sum_{i=1}^3 \deg \theta_i=|\A|$.

\item[(2)]
Let $H \in \A$. Then 
the derivations $\theta_2,\theta_3 \in D(\A^H,m^H)$ form a basis for $D(\A^H,m^H)$ if and only if 
$\det (\theta_i(x_j))=cQ(\A^H,m^H)$ up to nozero scalar $c \in \K$. In particular, 
$\sum_{i=2}^3 \deg \theta_i=|m|=:\sum_{X \in \A^H} m^H(X)=|\A|-1$.
\end{itemize}
\label{Saito}
\end{theorem}

\begin{theorem}[Terao's factorization, \cite{T2}, Main Theorem]
If $\A$ is free with $\exp(\A)=(1,a,b)$, then $$
\chi(\A;t)=(t-1)(t-a)(t-b).
$$
\label{Teraofactorization}
\end{theorem}

\begin{theorem}[\cite{A}, Theorem 1.1 and Corollary 1.2]
Assume that $\chi_0(\A;t)$ has two real roots $\alpha\le \beta$ and $H \in \A$. 
Then $|\A^H|-1 \le \alpha $ or $|\A^H|-1=\beta$. Moreover, $\A$ is free if $\alpha=|\A^H|-1$ or $\beta=|\A^H|-1$. 
\label{roots}
\end{theorem}

\begin{theorem}[\cite{Z}]
Assume that $\A$ is free with $\exp(\A)=(1,a,b)$. Then the Ziegler
restriction map $\pi^H:D_H(\A) \rightarrow D(\A^H,m^H)$ is surjective. 
In particular, $(\A^H,m^H)$ is free with $\exp(\A^H,m^H)=(a,b)$.
\label{Ziegler}
\end{theorem}

\begin{theorem}[Yoshinaga's criterion, \cite{Y2}, Theorem 3.2]
$$
\dim \coker \pi^H=\chi_0(\A;0)-d_1d_2,
$$
where $\exp(\A^H,m^H)=(d_1,d_2)$. Moreover, $\A$ is free with $\exp(\A)=(1,d_1,d_2)$ if and only if 
$\chi_0(\A;0)=d_1d_2$. 
\label{Ycriterion}
\end{theorem}


\begin{theorem}[Free surjection theorem, \cite{A12}, Theorem 3.3]
Let $H\in\A$.
If $\A'=\A\setminus\{H\}$ is free, then 
$$\rho^H:D(\A) \rightarrow D(\A^H)$$ is surjective.
\label{FST}
\end{theorem}


\begin{theorem}[\cite{A5}, Theorem 1.4]
Assume that $\A$ is free with exponents $(1,a,b)$ and $H \in \A$. If $\A':=
\A \setminus \{H\}$ is not free, then $\A'$ is SPOG with exponents $\POexp(\A)=(1,a,b)$ and level 
$d:=|\A'|-|\A^H|$, i.e., $D(\A')$ has the following minimal free resolution:
$$
0 \rightarrow S[-d-1] \rightarrow S[-1] \oplus S[-a] \oplus S[-b] \oplus S[-d] \rightarrow 
D(\A') \rightarrow 0.
$$
\label{SPOG}
\end{theorem}

Let us recall a part of a general theory of multiarrangements in $\K^2$. 
Let $\A$ be an arrangement in $W=\K^2$ and $m:\A \rightarrow \Z_{>0}$ be a multiplicity. 
For the coordinate ring $T=\mbox{Sym}^*(W^*)$, define 
$$
D(\A,m):=\{\theta \in \Der T \mid \theta(\alpha_H) \in T \alpha_H^{m(H)}\ (\forall H \in \A)\}.
$$
It is known, for example by \cite{Z}, that $D(\A,m)$ is reflexive, so free in this setup. Let $\exp(\A,m)=(a,b)$ be the set of degrees of a homogeneous free basis for $D(\A,m)$. Then we have the following.

\begin{lemma}[\cite{AN}, Lemma 4.2]
Let $\A$ be an arrangement in $\K^2$ and $m$ a multiplicity on $\A$. Let $H \in \A $ and 
$\delta_H:\A \rightarrow \{0,1\}$ be defined by $\delta_{H}(L)=1$ only when $L=H$ and $0$ otherwise. Let $\exp(\A,m-\delta_H)=(a,b)$. Then $\exp(\A,m)=(a+1,b)$ or $(a,b+1)$.
\label{AN}
\end{lemma}

\section{Proof of Theorem \ref{LPbound}}
In this section we prove Theorem \ref{LPbound}.
\medskip

\noindent
\textbf{Proof of Theorem \ref{LPbound}}.
(1)\,\, 
Let $\theta_E^H:=\rho^H(\theta_E)$ and let 
$$
D_1:=\overline{S} \theta_E^H + D(\A^H,m^H) \subset D(\A^H).
$$
Note that 
$$
\rho^H(D(\A))=\rho^H(S\theta_E \oplus D_H(\A)) \subset D_1 \subset D(\A^H).
$$
First let us consider $D_1/\rho^H(D(\A))$. By definition it is clear that 
$$
D_1/\rho^H(D(\A)) \simeq D(\A^H,m^H)/\rho^H(D(\A)) \cap D(\A^H,m^H).
$$
Since there is the Ziegler restriction map $\pi^H:D_H(\A) \rightarrow D(\A^H,m^H)$ and 
$\rho^H|_{D_H(\A)}=\pi^H|_{D_H(\A)}$, 
there is a surjection 
\begin{equation}
D(\A^H,m^H)/\rho^H(D_H(\A)) \rightarrow D(\A^H,m^H)/\rho^H(D(\A)) \cap D(\A^H,m^H).
\label{1iso}
\end{equation}By Theorem \ref{Ycriterion}, 
the left hand side is of dimension $\chi_0(\A;0)-d_1d_2$, where 
$(d_1,d_2)=\exp(\A^H,m^H)$. So
\begin{equation}
\chi_0(\A;0)-d_1d_2 \ge \dim D_1/\rho^H(D(\A)).
\label{Ybdd}
\end{equation}
Next let us consider $D(\A^H)/D_1$. 
We show that $\dim D(\A^H)/D_1=(d_1-n)(d_2-n)$. 
Let $\theta_E^H,\theta$ be a basis for $D(\A^H)$ with 
$\deg \theta=|\A^H|-1=:n$. Then for the basis $\theta_1,\theta_2$ for $D(\A^H,m^H)$ with $\deg 
\theta_i=d_i$, there is $a_i,b_i \in \overline{S}$ such that 
$$
\theta_i=a_i \theta_E^H+b_i \theta\ (i=1,2).
$$

First we consider the case $b_1=0$.  
Then $b_2 \neq 0$ by Theorem \ref{Saito}. Since $\theta_E^H$ is 
tangent to all hyperplanes with multiplicity one, it holds that $a_1=Q(\A^H,m^H)/Q(\A^H)$. 
Thus $d_1=|m^H|-|\A^H|+1$ and $d_2=|\A^H|-1=n$. So we may assume that 
$\theta=\theta_2$, i.e., $b_1=a_2=0$ and $b_2=1$. 
Hence
$$
D(\A^H)/D_1 \simeq 
\left({\overline{S}}\theta_E^H+{\overline{S}}\theta\right)/
\left({\overline{S}}\theta_E^H+{\overline{S}}a_1\theta_E^H+{\overline{S}}\theta\right)
=(0),
$$
and
\begin{equation}
\dim D(\A^H)/D_1=0=(d_1-n)(d_2-n)
\label{DD11}
\end{equation}
since $d_2=n$.

Next we consider the case $b_1b_2 \neq 0$.  Assume that 
there is a non-constant common divisor $b$ of $b_1$ and $b_2$. Since 
$Q(\A^H,m^H)=(a_1b_2-a_2b_1)Q(\A^H)$ up to non-zero scalar, we may assume that 
$b=\alpha_X$ for some $X \in \A^H$. Thus $m^H(X) \ge 2$. Since $\theta(\alpha_X) \in 
\overline{S}\alpha_X$ and $\theta_E^H(\alpha_X) =\alpha_X$, 
it holds that $\alpha_X \mid a_i$ for $i=1,2$. Thus $\alpha_X \mid \theta_i$ and 
$\theta_i/\alpha_X \in D(\A^H,m^H-\delta_X)$ for $i=1,2$. Since $\deg \theta_1/\alpha_X+
\deg \theta_2/\alpha_X=|m^H|-2<|m^H-\delta_X|$, and $\theta_1/\alpha_X,\ 
\theta_2/\alpha_X$ are $\overline{S}$-independent, this contradicts to Theorem \ref{Saito}. 
Thus $(b_1,b_2)=1$. 
Hence 
$$
D(\A^H)/D_1 \simeq 
\left({\overline{S}}\theta_E^H+{\overline{S}}\theta\right)/
\left({\overline{S}}\theta_E^H+{\overline{S}}b_1\theta+{\overline{S}}b_2\theta\right)
\simeq \overline{S}/(b_1,b_2)
.
$$
Since $\deg b_i=d_i-n$ for $i=1,2$, 
they are complete intersection. Thus 
\begin{equation}
\dim D(\A^H)/D_1=(d_1-n)(d_2-n). 
\label{DD12}
\end{equation}

Now using (\ref{Ybdd}), (\ref{DD11}) and (\ref{DD12}), 
\begin{eqnarray}\label{eqsA}
\dim D(\A^H)/\rho^H(D(\A))&=&
\dim D(\A^H)/D_1 +\dim D_1/\rho^H(D(\A))\nonumber\\
&\le&(n-d_1)(n-d_2)+\chi_0(\A;0)-d_1d_2\\
&=&n^2-(d_1+d_2)n+\chi_0(\A;0)=\chi_0(\A;n)\nonumber,
\end{eqnarray}
which completes the proof. 

(2)\,\,
If $\A$ is free, 
we know that $\rho^H(D_H(\A))=D(\A^H,m^H)$
and 
$(a,b)=(d_1,d_2)=\exp(\A^H,m^H)$ by Theorem \ref{Ziegler}. 
Thus in (\ref{1iso}), 
we have a surjection 
$$
0 \rightarrow D_1/\rho^H(D(\A))=0.
$$
So the inequalities (\ref{Ybdd}) and that in 
(\ref{eqsA}) are the equalities.
By the argument in the proof of (1), we know that for 
$$
\coker \rho^H=D(\A^H)/\rho^H(D(\A))=
D(\A^H)/D_1 
\simeq \overline{S}/(b_1,b_2)[-n].
$$
Thus $(b_1,b_2)$ form a regular sequence. Thus 
\begin{eqnarray*}
0 \rightarrow \overline{S}[-a-b+2n] \rightarrow \overline{S}[-a+n]\oplus 
\overline{S}[-b+n ]
\rightarrow \overline{S} \rightarrow \coker \rho^H \rightarrow 0
\end{eqnarray*}
is a minimal free resolution of $\coker \rho^H$, and we have 
$$
\hilb(
\coker\rho^H
;x)=x^n \displaystyle \frac{1-x^{a-n}}{1-x}
\displaystyle \frac{1-x^{b-n}}{1-x}.
\eqno{\square}
$$
\medskip

Based on Theorem \ref{LPbound}, we can show the higher version. For that recall the following:

\begin{define}
Let $\ell=2$ or $3$ and $\A$ be an arrangement in $\K^\ell$. Then 
\begin{itemize}
    \item [(1)]
$D^2(\A)$ is defined by 
$$
D^2(\A):=\{\theta \in \wedge^2 \Der S \mid \theta(\alpha_H,f) \in S\alpha_H \ (\forall H \in \A)\}.
$$
\item [(2)] 
The Euler restriction map $\rho_2^H:D^2(\A) \rightarrow D^2(\A^H)$ is defined by 
$$
\rho_2^H(\theta)(\overline{f},\overline{g}):=\overline{{\theta(f,g)}}.
$$
\end{itemize}
\end{define}

On $D^2(\A)$, the following is known.

\begin{theorem}[\cite{ST}, Proposition 3.4]
Let $\A$ be an arrangement in $\K^\ell$ with $\ell \in \{2,3\}$. 
\begin{itemize}

    \item[(1)]
    Assume that $\A$ is free. Then $D^2(\A) =\wedge^2 D(\A)$.
    \item[(2)]
    Let $\ell=2$. Then $D^2(\A)$ is a free $S$-module with a basis $Q(\A)\partial_{x_1} \wedge \partial_{x_2}$.
\end{itemize}
\end{theorem}

Now we prove the same result as in Theorem \ref{LPbound} for $\rho_2^H$.

\begin{theorem}
\begin{itemize}
    \item[(1)]
$$
\dim \coker \rho_2^H \le LP(\A;H),
$$
and 
\item[(2)]
$$
\dim \coker \rho^H_2 = LP(\A;H),
$$
if $\A$ is free. Moreover, 
    $$
    \hilb(\coker \rho_2^H;x)=x^{|\A^H|-1}
    \displaystyle \frac{1-x^{a-|\A^H|+1}}{1-x}
    \displaystyle \frac{1-x^{b-|\A^H|+1}}{1-x},
    $$
    where $\exp(\A)=(1,a,b)$.
\end{itemize}
\label{LPboud2}
\end{theorem}

\noindent
\textbf{Proof}. 
(1)\,\,
Let $\alpha_H=x_1$ and we may assume that $x_2=0$ is in $\A^H$.
Then $\theta_E^H, \theta_2:=(Q(\A^H)/x_2)\partial_{x_3}$ form a basis 
for $D(\A^H)$. Also, $\theta_E^H \wedge \theta_2=Q(\A^H) \partial_{x_2} \wedge \partial_{x_3}$ form a basis for $D^2(\A^H)$. Let $f_1 \theta_E^H+f_2 \theta_2 \in 
\rho_2^H(D(\A))$. Then there is $\varphi \in D(\A)$ such that  $\rho_2^H(\varphi)=f_2\theta_2$. Then 
$$
(f_1 \theta_E^H+f_2 \theta_2) \wedge \theta_E^H=f_2 
Q(\A^H) \partial_{x_2} \wedge \partial_{x_3} 
=\rho_2^H(\theta_E \wedge \varphi) \in \rho_2^H(D^2(\A^H)).
$$
Hence there is a well-defined surjection 
$$\wedge \theta_E^H\colon D(\A^H)/\rho_2^H(D(\A)) \rightarrow D^2(\A^H)/\rho_2^H(D^2(\A))$$ 
by taking the wedge with $\theta_E^H$. 
Consider the kernel of the map $\wedge \theta_E^H$. Let $\theta \in D(\A^H)$ satisfy that $\theta \wedge \theta_E^H \in \rho_2^H(D^2(\A))$. Also, for $L \in \A^H$, 
define the map $\partial_L:D^2(\A^H) \rightarrow D(\A^H)$ by 
$$
\partial_L(\theta)(f):=\theta(\alpha_L,f)/\alpha_L
$$
for $\theta \in D^2(\A^H),\ f \in \overline{S}$. 
Since $$
\partial_L (\rho_2^H(D^2(\A))) \ni \partial_L(\rho_2^H(\theta))=\rho_2^H(\theta)(\alpha_L,*)/\alpha_L=
\rho_2^H(\partial_L(\theta)) \in \rho_2^H(D(\A))
$$
for $\theta \in D^2(\A)$, we know that 
$$
\partial_L:D^2(\A^H)/\rho_2^H(D^2(\A)) \rightarrow D(\A^H)/\rho_2^H(D(\A))
$$
is well-defined. Let $f_1 \theta_E^H+f_2\theta_2 \in D(\A^H)/\rho_2^H(D(\A))$. 
We may assume that $\theta_2 \in D_L(\A^H)$. 
Then $$
\partial_L (\theta_E^H \wedge (f_1 \theta_E^H+f_2\theta_2))=
\pm f_2 \theta_2,
$$
which is the same as $f_1 \theta_E^H+f_2\theta_2$ since $
\theta_E^H \in \rho^H_2(D(\A^H))$. Hence $\wedge \theta_E^H$ is injective, and 
we know that 
$$
D(\A^H)/\rho_2^H(D(\A)) \simeq D^2(\A^H)/\rho_2^H(D^2(\A)).
$$
Combining this with Theorem \ref{LPbound}, we know that $\dim \coker \rho_2^H \le LP(\A;H)$. 

(2)\,\,
The same proof as in Theorem \ref{LPbound} and the argument in (1) complete the proof. \owari
\medskip

Theorem \ref{LPbound} (2) is not true in general for non-free arrangements. 

\begin{example}
Let $\A$ be $$
xyz(x+y+z)=0.
$$
Then $|\A^H|=3$ for any $H \in \A$ and 
$\chi_0(\A;t)=t^2-3t+3$. Hence 
$$
\chi_0(\A;|\A^H|-1)=1.
$$
However, it is easy to show that $\rho^H$ is surjective for all $H$. Thus 
$\dim \coker \rho^H=0<1=\chi_0(\A;2)$.

\label{generic}
\end{example}

\section{Multiarrangements and characteristic polynomials}

In the previous section we showed the importance of the combinatorial invariant $\chi_0(\A;|\A^H|-1)$. Not only $\coker (\rho^H:D(\A) \rightarrow D(\A^H))$, but also 
$$
\coker (\pi^H\colon D_H(\A) \rightarrow D(\A^H,m^H)),
$$
which is explicitly written by Yoshinaga's criterion, 
is bounded by LP-invariants as follows:

\begin{theorem}
\begin{eqnarray*}
\chi_0(\A;|\A^H|-1)&=&\chi_0(\A;|\A|-|\A^H|)
 \ge 
 \dim \coker (\pi^H:D_H(\A) \rightarrow D(\A^H,m^H))\\
 &=&\chi_0(\A;0)-d_1d_2.
\end{eqnarray*}
where $\exp(\A^H,m^H)=(d_1,d_2)$.
\label{Zieglerdim}
\end{theorem}

\noindent
\textbf{Proof}. 
Compute 
\begin{eqnarray*}
\chi_0(\A;|\A|-|\A^H|)&=&
(|\A|-|\A^H|)^2-(|\A|-1)(|\A|-|\A^H|)+b_2^0\\
&=&
b_2^0-(|\A^H|-1)(|\A|-|\A^H|).
\end{eqnarray*}
Note that $\exp(\A^H)=(1,|\A^H|-1)$. By Lemma \ref{AN}, we know that 
if $\exp(\B,m)=(a,b)$ for some $2$-arrangement $\B$ and a multiplicity $m$ on it, 
then for $H \in \B$, it holds that $\exp(\B,m+\delta_H)=(a+1,b)$ or $(a,b+1)$. Since 
$\exp(\A^H,k)=(k,|\A^H|-1)$ if $|k| \le 2|\A^H|-1$, the minimum of $d_1d_2$ is attained when $d_1=|\A^H|-1$ and $d_2=|\A|-|\A^H|$. As a conclusion,
\begin{eqnarray*}
\chi_0(\A;|\A|-|\A^H|)&=&
b_2^0-(|\A^H|-1)(|\A|-|\A^H|)
\ge b_2^0-d_1d_2.
\end{eqnarray*}
So Yoshinaga's criterion completes the proof. \owari
\medskip

As we saw, the dimension of the cokernel of the 
Ziegler restriction is explicitly computable by Yoshinaga's criterion, but for that we need to know the exponents of multi-arrangements, which is in general hard to see. Theorem \ref{Zieglerdim} gives us an easy way to estimate 
it just based on the combinatorial data. Moreover, 
we can understand Yoshinaga's criterion as a ``value'' of the reduced characteristic polynomial with some 
integer substituted:

\begin{prop}
Let $\exp(\A^H,m^H)=(d_1,d_2)$. Then 
$$
\dim \coker \pi^H=\chi_0(\A;d_1)=\chi_0(\A;d_2).
$$
\label{Ychar}
\end{prop}

\noindent
\textbf{Proof}.
Compute 
\begin{eqnarray*}
\chi_0(\A;d_1)&=&
d_1^2-(|\A|-1)d_1+b_2^0
=b_2^0-d_1d_2.
\end{eqnarray*}
So Yoshinaga's criterion completes the proof. \owari
\medskip

There is no obvious relation between $\dim \ker \rho^H$ and $\dim \ker \pi^H$ as 
seen below.

\begin{example}
(1)\,\,
First consider the case when 
$\dim \coker \rho^H>\dim \coker \pi^H$. 
Let $\A$ be the cone of the affine arrangement in $\R^2$ consisting of the edges and diagonals of the regular pentagon, which is thus an arrangement in $\R^3$. It is known to be free with exponents $(1,5,5)$, and 
$|\A^H|=5$ for all $H \in \A$. Thus $LP(\A,H)=(4-5)^2=1$. Since it is free, $\dim \coker \pi^H=0$ by Theorem \ref{Ycriterion}. Since $\A$ is free, Theorem \ref{LPbound} shows that $\dim \coker \rho^H=LP(\A,H)=1>
\dim \coker \pi^H=0$.

(2)\,\,
Second consider the case when 
$\dim \coker \rho^H<\dim \coker \pi^H$. 
Let $\A$ be defined by $xyz(x+y+z)=0$ in $\R^3$.
Then $\chi_0(\A;t)=t^2-3t+3$ thus $\A$ is not free. Since $|\A^H|=3$ for all $H \in \A$, it holds that $LP(\A,H)=1$. Since $\A$ is not free, Theorem \ref{Zieglerdim} shows that $\dim \coker \pi^H=1$. On the other hand, 
it is easy to compute that $\dim \coker \rho^H=0<1=\dim \coker \pi^H$.

(3)\,\,
Third consider the case when 
$\dim \coker \rho^H=LP(\A,H)$ but $\A$ is not free. 
Let $\A$ be defined by $$
xyz(x-z)(y-z)(x-y+z)(x-y-z)=0
$$
in $\R^3$. Then $\chi_0(\A;t)=t^2-6t+10$, thus $\A$ is not free and 
$\dim \coker \pi^H>0 $ by Theorem \ref{Ycriterion}. Since we can compute that $\exp(\A^H,m^H)=(3,3)$, Theorem \ref{Ycriterion} shows that $\dim \coker \pi^H =10-9=1$.
Let $H:z=0$. Then $|\A^H|=3$ and thus 
$LP(\A,H)=2$. Since $\A \cup \{x-y=0\}$ is free with exponents 
$(1,3,4)$, Theorem \ref{SPOG} shows that $D(\A)$ has a minimal free resolution of the following form:
$$
0 \rightarrow S[-5] 
\rightarrow S[-1] \oplus S[-3] \oplus S[-4]^2 \rightarrow D(\A) 
\rightarrow 0.
$$
Since $\exp(\A^H)=(1,2)$, we can compute that $2\le \dim \coker \rho^H$. By Theorem \ref{LPbound}, it holds that 
$$
LP(\A,H)=\dim \coker \rho^H =2>1=\dim \coker \pi^H .
$$

\end{example}

By using these, we can give an application to determine the structure of $D(\A)$ combinatorially.

\begin{cor}
Let $\ell=3$ and assume that $LP(\A,H)=1$ for some $H \in \A$. Then either $\A$ is free, or nearly free with exponents $(a,b)$. Namely, there is a minimal free resolution 
\begin{equation}
0 \rightarrow S[-b-1] 
\rightarrow S[-1] \oplus S[-a] \oplus S[-b]^2 \rightarrow 
D(\A) \rightarrow 0.
\label{eq1}
\end{equation}
\label{cor1}
\end{cor}

\begin{example}
Let $\A$ be an arrangement in $\R^3$ defined by $$
xz(y-z)(x^2-y^2)(x^2-4y^2)(x-3y)=0.
$$
Then $\chi_0(\A;t)=t^2-7t+11$. By Theorem \ref{Teraofactorization}, $\A$ is not free. Let $H:x=0$. Then $|\A^H|=3$ and hence $L(\A,H)=\chi_0(\A;2)=1$. Hence Corollary \ref{cor1} shows that 
$D(\A)$ has a minimal free resolution of the type (\ref{eq1}). In this case, explicitly, $a=2$ and $b=6$.
\end{example}

%
%

%
%
%


\begin{thebibliography}{ABCHT}


\bibitem{A}
T. Abe, 
Roots of characteristic polynomials and 
and intersection points of line arrangements. 
\textit{J. Singularities}, 
\textbf{8} (2014), 100--117.







\bibitem{A5}
T. Abe, 
Plus-one generated and next to free arrangements of hyperplanes. 
\textit{Int. Math. Res. Not.} \textbf{2021} (2021), no. 12, 9233--9261.



\bibitem{A8}
T. Abe, 
Double points of free projective line arrangements, 
\textit{Int. Math. Res. Not}.
\textbf{2022} (2022), no. 3, 1811--1824.






\bibitem{A12}
T. Abe, 
Addition-deletion theorems for the Solomon-Terao polynomials and 
$B$-sequences of hyperplane arrangements. 
\textit{Math. Z.} \textbf{306}(2024), no. 2, 25.























\bibitem{AN}
T. Abe and Y. Numata,
Exponents of 2-multiarrangements and multiplicity lattices. 
J. Alg. Combin. \textbf{35} (2012), 1--17.



















\bibitem{OT} P. Orlik and H. Terao, \textit{Arrangements of hyperplanes}.
Grundlehren der Mathematischen Wissenschaften, 
\textbf{300}. Springer-Verlag, Berlin, 1992.



\bibitem{Sa}
K. Saito, 
Theory of logarithmic differential forms and logarithmic vector fields.
\textit{J. Fac. Sci. Univ. Tokyo} \textbf{27} (1980), 265--291.   





\bibitem{ST}
L. Solomon and H. Terao, 
A formula for the characteristic polynomial
of an arrangement. \textit{Adv. Math.} \textbf{64} (1987), 
305--325.


\bibitem{T2}
H. Terao, 
Generalized exponents of a free arrangement of hyperplanes and
Shephard-Todd-Brieskorn formula. \textit{Invent. math}. 
\textbf{63}  (1981),
159--179.



\bibitem{Y2}
M. Yoshinaga, 
On the freeness of 3-arrangements. 
\textit{Bull. London Math. Soc.} \textbf{37} (2005), no. 1, 126--134. 



\bibitem{Z}
G. M. Ziegler, 
Multiarrangements of hyperplanes and their freeness.  Singularities (Iowa City, IA, 1986),  345--359,
Contemp. Math., {\bf 90}, Amer. Math. Soc., Providence, RI, 1989. 


	
\end{thebibliography}
\end{document}